\magnification=\magstep1
\input amstex
\NoBlackBoxes
\documentstyle{amsppt}

\define\C{{\Bbb C}}
\define\R{{\Bbb R}}
\define\dee{\partial}

\redefine\O{\Omega}
\redefine\phi{\varphi}
\define\Obar{\overline{\Omega}}
\define\Ot{\widetilde\Omega}

\NoRunningHeads
\topmatter
\title
Complexity of the classical kernel functions of potential theory
\endtitle
\author Steven R. Bell${}^*$ \endauthor
\thanks ${}^*$Research supported by NSF grant DMS-9302513 \endthanks
\keywords Bergman kernel, Szeg\H o kernel, Green's funcion, Poisson kernel
\endkeywords
\subjclass 32H10 \endsubjclass
\address
Mathematics Department, Purdue University, West Lafayette, IN  47907 USA
\endaddress
\email bell\@math.purdue.edu \endemail
\abstract
We show that the Bergman, Szeg\H o, and Poisson kernels associated to a
finitely connected domain in the plane are all composed of finitely many
easily computed {\it functions of one variable}.  The new formulas give
rise to new methods for computing the Bergman and Szeg\H o kernels in which
all integrals used in the computations are line integrals; at no point is an
integral with respect to area measure required.  The results mentioned so
far can be interpreted as saying that the kernel functions are simpler than one
might expect.  However, we also prove that the kernels cannot be too simple
by showing that the only finitely connected domains in the plane whose
Bergman or Szeg\H o kernels are rational functions are the obvious ones.  This
leads to a proof that the classical Green's function associated to a finitely
connected domain in the plane is the logarithm of a rational function if and
only if the domain is simply connected and rationally equivalent to the unit
disc.
\endabstract
\endtopmatter
\document

\hyphenation{bi-hol-o-mor-phic}
\hyphenation{hol-o-mor-phic}

\subhead 1. Introduction \endsubhead
The Bergman and Szeg\H o kernels associated to a bounded domain in
the plane with smooth boundary carry encoded within them an
astonishing amount of information about the domain.  Conformal
mappings onto canonical domains, classical domain functions, and
other important objects of potential theory can be expressed simply
in terms of the Bergman and Szeg\H o kernels.  It is therefore
tempting to believe that these kernels are extremely complex and
difficult to compute.  The purpose of this paper is to show that
the kernel functions are not nearly as complex as one might suspect.

Suppose that $\O$ is a bounded finitely connected domain in the
plane with $C^\infty$ smooth boundary, i.e., that the boundary $b\O$
of $\O$ is given by finitely many non-intersecting $C^\infty$ simple
closed curves.  The Bergman kernel $K(z,w)$ and the Szeg\H o kernel
$S(z,w)$ associated to such a domain are both known to extend to be in
the space $C^\infty((\Obar\times\Obar)-\Cal D)$ where $\Cal D$
denotes the boundary diagonal $\{(z,z):z\in b\O\}$.  Our problem is
to determine a method to compute $K(z,w)$ and $S(z,w)$ at any given
ordered pair of points $(z,w)$ in $(\Obar\times\Obar)-\Cal D$.  We shall
see that, once the boundary values of finitely many basic functions of one
variable have been determined, the kernels become known at all points
$(z,w)$.  Furthermore, the basic functions which comprise the
kernel functions are all solutions to explicit Kerzman-Stein integral
equations, and as such, are easy to compute.  All elements
of the kernel functions may be computed by means of simple linear
algebra and one dimensional integrals and one dimensional integral equations.
At no point is a double integral needed.

Although the words, ``numerical method,'' appear in this paper, this is
not a paper on numerical analysis; no examples of numerical computations
are given.  However, the results of this paper should be interesting to
numerical analysts.

Our results are particularly surprising for the Bergman kernel function.
A traditional way to attempt to compute the Bergman kernel has been to
orthonormalize a set of rational functions that span a dense subset of
the Bergman space.  This is a numerical nightmare compared to the methods
we establish below.

In \S6, we give conditions on a domain
for its Bergman or Szeg\H o kernel function to be a rational function.
We prove, for example, that the Bergman kernel associated to a
finitely connected domain is rational exactly when the domain
is simply connected and is biholomorphic to the unit disc via a rational
mapping.  Thus, the domains whose Bergman kernels are rational are precisely
those domains whose kernels can be seen to be rational by means of an
elementary application of the transformation formula for the kernel functions
under biholomorphic maps and the fact that the kernel for the unit
disc is rational.  This result about the Bergman kernel has as a
corollary that the Green's function associated to a finitely
connected domain is the logarithm of a rational function if and only
if the domain is simply connected and is rationally equivalent to the
unit disc.

In \S7, we show how the Poisson kernel can be expressed in terms of the
Szeg\H o kernel, and thereby shed light on the degree of complexity of
the Poisson kernel.

Our results are most interesting in case the domain under study is multiply
connected.  However, to illustrate the flavor of our results, we take
a moment here to state analogues of our theorems for a bounded
simply connected domain $\O$ with $C^\infty$ smooth boundary.  How
difficult is it to compute the Szeg\H o kernel $S(z,w)$ at any given
pair of points?  Is it so difficult that we would need to follow a
separate numerical procedure to compute $h(z)=S(z,w_0)$ for each
individual point $w_0\in\O$?  The following formula shows that it is not
nearly so difficult.  It is only as difficult as computing the boundary
values of a single function.  Let $a$ be a fixed point in $\O$ and let
$f_a(z)$ denote the Riemann mapping function mapping $\O$ one-to-one onto
the unit disc $D_1(0)$ with $f_a(a)=0$ and $f_a'(a)>0$.  This Riemann map can
easily be written down in terms of the boundary values of the function
$S(z,a)$ (see \cite{2}), and $S(z,a)$ is the solution to a simple
Fredholm integral equation of the second kind with $C^\infty$ kernel
and inhomogeneous term (see \cite{{2,3,7,10,11,14}}).  The kernel $S(z,w)$
may be expressed as
$$S(z,w)=\frac{c\,S(z,a)\overline{S(w,a)}}{1-f_a(z)\overline{f_a(w)}},$$
where $c=1/S(a,a)$.  This shows that, once the boundary values of the
single function of one variable $S(z,a)$ are known, the Szeg\H o kernel
can be evaluated at an arbitrary pair of points.  A similar identity holds
for the Bergman kernel,
$$K(z,w)=\frac{4\pi S(z,a)\overline{S(w,a)}}{(1-f_a(z)\overline{f_a(w)})^2}.$$
This shows that the Bergman kernel is composed of the same basic
functions that make up the Szeg\H o kernel.  Finally, the Poisson
kernel $p(z,w)$ is given by
$$p(z,w)=
\frac{S(z,w)S(w,a)}{S(z,a)}+
\frac{\overline{S(z,w)S(w,a)f_a(z)}}{\overline{S(z,a)f_a(w)}},$$
where $z$ is a point in $\O$ and $w$ is a point in the boundary
(see \cite{2, page~37}).  Thus, the Poisson kernel is also composed of the
same basic functions.  None of these formulas for the kernel functions in
a simply connected domain could be considered very new.  However, we shall
prove analogous results for $n$-connected domains that are new.  The new
results show that there are $n+1$ basic functions that comprise all the
kernels.  An interesting feature of all the results in this paper is the
central role played by the zeroes of the Szeg\H o kernel.

\subhead 2. The Ahlfors map and zeroes of the Szeg\H o kernel \endsubhead
Before we start stating and proving our main theorems, we must review
some basic facts about the kernel functions.  Most of these facts
are proved in Bergman's book \cite{5}; all of them are proved
in \cite{2}.

Suppose that $\O$ is a bounded $n$-connected domain in the
plane with $C^\infty$ smooth boundary.  Let $\gamma_j$, $j=1,\dots,n$,
denote the $n$ non-intersecting $C^\infty$ simple closed curves which
define the boundary of $\O$, and suppose that $\gamma_j$ is
parameterized in the standard sense by $z_j(t)$, $0\le t\le 1$.  We
shall use the convention that $\gamma_n$ denotes the {\it outer
boundary curve\/} of $\O$.  Let $T(z)$ be the $C^\infty$ function
defined on $b\O$ such that $T(z)$ is the complex number representing
the unit tangent vector at $z\in b\O$ pointing in the direction of
the standard orientation.  This complex unit tangent vector function
is characterized by the equation $T(z_j(t))=z_j'(t)/|z_j'(t)|$.

We shall let $A^\infty(\O)$ denote the space of holomorphic functions
on $\O$ that are in $C^\infty(\Obar)$.  The space of complex valued
functions on $\O$ that are square integrable with respect to Lebesgue
area measure $dA$  will be written $L^2(\O)$, and the space of complex
valued functions on $b\O$ that are square integrable with respect to arc
length measure $ds$ will be denoted by $L^2(b\O)$.  The Bergman space of
holomorphic functions on $\O$ that are in $L^2(\O)$ shall be written
$H^2(\O)$ and the Hardy space of functions in $L^2(b\O)$ that are the
$L^2$ boundary values of holomorphic functions on $\O$ shall be written
$H^2(b\O)$.  The inner products associated to $L^2(\O)$ and $L^2(b\O)$
shall be written
$$\langle u,v\rangle_{\O}=\iint_{\O}u\ \bar v\ dA\quad\text{ and }\quad
\langle u,v\rangle_{b\O}=\int_{b\O}u\ \bar v\ ds,$$
respectively.

For each fixed point $a\in\O$, the Szeg\H o kernel $S(z,a)$, as a function of
$z$, extends to the boundary to be a function in $A^\infty(\O)$.  (An even
stronger smoothness property is mentioned in the introduction.) \ Furthermore,
$S(z,a)$ has exactly $(n-1)$ zeroes in $\O$ (counting multiplicities) and
does not vanish at any points $z$ in the boundary of $\O$.  The {\it
Garabedian kernel\/} $L(z,a)$ is a kernel related to the Szeg\H o
kernel via the identity
$$\frac{1}{i} L(z,a)T(z)=S(a,z)\qquad\text{for $z\in b\O$ and $a\in\O$.}
\tag2.1$$
For fixed $a\in\O$, the kernel $L(z,a)$ is a holomorphic function of $z$
on $\O-\{a\}$ with a simple pole at $a$ with residue $1/(2\pi)$. 
Furthermore, as a function of $z$, $L(z,a)$ extends to the boundary
and is in the space $C^\infty(\Obar-\{a\})$.  In fact, $L(z,a)$
extends to be in $C^\infty((\Obar\times\Obar)-\{(z,z):z\in\Obar\})$.  Also,
$L(z,a)$ is non-zero for all $(z,a)$ in $\Obar\times\O$ with $z\ne a$.

The kernel $S(z,w)$ is holomorphic in $z$ and antiholomorphic in $w$
on $\O\times\O$, and $L(z,w)$ is holomorphic in both variables for
$z,w\in\O$, $z\ne w$.  We shall need to know that $S(z,z)$ is real and
positive for each $z\in\O$, and we shall need to use the basic
identities $S(z,w)=\overline{S(w,z)}$ and $L(z,w)=-L(w,z)$.  The
Szeg\H o kernel reproduces holomorphic functions in the sense that
$$h(a)=\langle h, S(\cdot,a)\rangle_{b\O}$$
for all $h\in H^2(b\O)$ and $a\in\O$.

Given a point $a\in\O$, the Ahlfors map $f_a$ associated to the pair ($\O,a)$
is a proper holomorphic mapping of $\O$ onto the unit disc.  It is an
$n$-to-one mapping (counting multiplicities), it extends to be in
$A^\infty(\O)$, and it maps each boundary curve $\gamma_j$ one-to-one
onto the unit circle.  Furthermore, $f_a(a)=0$, and $f_a$ is the unique
function mapping $\O$ into the unit disc maximizing the quantity $|f_a'(a)|$
with $f_a'(a)>0$.  The Ahlfors map is related to the Szeg\H o kernel
and Garabedian kernel via
$$f_a(z)=\frac{S(z,a)}{L(z,a)}.\tag2.2$$
Also, $f_a'(a)=2\pi S(a,a)\ne 0$.  Because $f_a$ has $n$ zeroes, and
because the simple pole of $L(z,a)$ at $a$ accounts for the simple zero of
$f_a$ at $a$, it follows that $S(z,a)$ has $(n-1)$ zeroes in
$\O-\{a\}$.  Let $a_1,a_2,\dots,a_{n-1}$ denote these $n-1$ zeroes
(counted with multiplicity).  I proved in \cite{4} (see also \cite{2,
page~105}) that, if $a$ is close to one of the boundary curves, then the zeroes
$a_1,\dots,a_{n-1}$ become distinct simple zeroes.  It follows from
this result that, for all but at most finitely many points $a\in\O$,
$S(z,a)$ has $n-1$ distinct simple zeroes in $\O$ as a function of $z$.

\subhead 3. A special orthonormal basis for the Hardy space \endsubhead
The zeroes of the Szeg\H o kernel give rise to a particularly nice
basis for the Hardy space of an $n$-connected domain with $C^\infty$
smooth boundary.  We shall use the notation that we set up in the
preceding section.  We assume that $a\in\O$ is a fixed point in $\O$
that has been chosen so that the $n-1$ zeroes, $a_1,\dots,a_{n-1}$,
of $S(z,a)$ are distinct and simple.  We shall let $a_0$ denote $a$
and we shall use the shorthand notation $f(z)$ for the Ahlfors map
$f_a(z)$.

We shall now prove that the set of functions
$\{h_{ik}(z): 0\le i\le n-1,\text{ and }k\ge 0\}$
where $h_{ik}$ is defined via
$$h_{ik}(z)=S(z,a_i)f(z)^k$$
forms a basis for the Hardy space $H^2(b\O)$.  Furthermore,
$$\langle h_{ik},h_{jm}\rangle_{b\O}=
\cases
0, &\text{if $k\ne m$} \\
S(a_j,a_i), &\text{if $k=m$}.\endcases\tag3.1$$
The proof of these assertions consists of two parts.  First,
we must prove that these functions span a dense subset of
$H^2(b\O)$, and second, we must prove the identity (3.1) to
see that the set forms a basis.  To prove the density of the span, suppose
that $g\in H^2(b\O)$ is orthogonal to the span.  Notice that
$$\langle g,S(\cdot,a_j)\rangle_{b\O}=g(a_j),$$
and therefore $g$ vanishes at $a_0,a_1,\dots,a_{n-1}$.
Suppose we have shown that $g$ vanishes to order $m$ at
each $a_j$, $j=0,1,\dots,n-1$.  It follows that $g/f^m$ is in
$H^2(b\O)$ and the value of $g/f^m$ at $a_j$ is $g^{(m)}(a_j)/f'(a_j)^m$.
Since $|f(z)|=1$ when $z\in b\O$, it follows that $1/f(z)=\overline{f(z)}$
when $z\in b\O$, and we may write
$$\langle g,S(\cdot,a_j)f^m\rangle_{b\O}=
\langle g/f^m,S(\cdot,a_j)\rangle_{b\O}=
g^{(m)}(a_j)/f'(a_j)^m.$$
We conclude that $g$ vanishes to order $m+1$ at each $a_j$.  By
induction, $g$ vanishes to infinite order at each $a_j$ and hence,
$g\equiv0$.  This proves the density.  To prove (3.1), let us suppose
first that $k>m$.  The fact that $\overline{f}=1/f$ on $b\O$ and
the reproducing property of the Szeg\H o kernel now yield that
$$\gather
\langle h_{ik},h_{jm}\rangle_{b\O}=
\int_{z\in b\O} S(z,a_i)f(z)^{k-m}\ \overline{S(z,a_j)}\ ds= \\
\int_{z\in b\O} S(a_j,z)\ S(z,a_i)f(z)^{k-m}\ ds=S(a_j,a_i)f(a_j)^{k-m}.
\endgather$$
The identity now follows because $f(a_j)=0$ for all $j$.  If $k=m$,
then
$$\langle h_{ik},h_{jm}\rangle_{b\O}=
\int_{z\in b\O} S(a_j,z)\ S(z,a_i)\ ds=S(a_j,a_i),$$
and identity (3.1) is proved.  We remark here that it easy to see that
the functions $h_{ik}$ are linearly independent.  Indeed, identity
(3.1) reveals that we need only check that, for fixed $k$,  the $n$
functions $h_{ik}$, $i=0,1,\dots,n-1$, are independent, and this is
easy because a relation of the form
$$\sum_{i=0}^{n-1} C_i S(z,a_i)\equiv 0$$
implies, via the reproducing property of the Szeg\H o kernel, that every
function $g$ in the Hardy space satisfies
$$\sum_{i=0}^{n-1} \overline{C_i}\,g(a_i)=0,$$
and it is easy to construct polynomials $g$ that violate such a
condition.

We next orthonormalize the sequence $\{h_{ik}\}$ via the Gram-Schmidt
procedure.  Identity (3.1) shows that most of the functions in the sequence
are already orthogonal, and so our task is quite easy.  We need only fix
$k$ and orthonormalize the $n$~functions $h_{ik}$, $i=0,1,\dots,n-1$.  We
obtain an orthonormal set $\{H_{ik}\}$ given by
$$\gather
H_{0k}(z)=b_{00} S(z,a)f(z)^k\qquad\text{and,} \\
H_{ik}(z)=\sum_{j=1}^{i}b_{ij} S(z,a_j)f(z)^k,\qquad i=1,\dots,n-1
\endgather$$
where $b_{ii}\ne0$ for each $i=0,1,\dots,n-1$.  Because $|f|=1$ on $b\O$,
it follows that {\it the coefficients\/} $b_{ij}$ {\it do not depend on\/}
$k$.  Notice that $H_{ij}$ does not contain a term involving $S(z,a)$ if $i>0$
because of (3.1) and the fact that $S(a_i,a)=0$.

The Szeg\H o kernel can be written in terms of our orthonormal basis
via
$$S(z,w)=\sum_{i=0}^{n-1}\sum_{k=0}^\infty H_{ik}(z)\,\overline{H_{ik}(w)}.$$
The sum
$$\sum_{k=0}^\infty f(z)^k\,\overline{f(w)^k}=\frac{1}{1-f(z)\overline{f(w)}}$$
can be factored from the expression for $S(z,w)$ to yield the formula in
the following theorem.
\proclaim{Theorem 3.1}
The Szeg\H o kernel can be evaluated at an arbitrary pair of points
$(z,w)$ in $\O\time\O$ via the formula
$$S(z,w)=\frac{1}{1-f(z)\overline{f(w)}}\left(c_0 S(z,a)\overline{S(w,a)}+
\sum_{i,j=1}^{n-1} c_{ij}S(z,a_i)\,\overline{S(w,a_j)}\right)\tag3.2$$
where $f(z)$ denotes the Ahlfors map $f_a(z)$, $c_0=1/S(a,a)$, and the
coefficients $c_{ij}$ are given as the coefficients of the inverse matrix
to the matrix $\left[S(a_j,a_k)\right]$.
\endproclaim

The only part of Theorem~3.1 that remains unproved is the statement
about the coefficients in the formula.  We have shown that these coefficients
exist and that they are given as certain combinations of the Gram-Schmidt
coefficients used above.  That $c_0=1/S(a,a)$ can be seen by setting $z=a$
and $w=a$ in (3.2).  To complete the proof of Theorem~3.1, we shall
now describe how to determine the coefficients $c_{ij}$.  Suppose
$1\le k\le n-1$.  Set $w=a_k$ in (3.2) and note that $f(a_k)=0$ and
$S(a,a_k)=0$ to obtain
$$S(z,a_k)=\sum_{i=1}^{n-1}
\left(\sum_{j=1}^{n-1}c_{ij}S(a_j,a_k)\right)S(z,a_i).$$
We saw an identity like this when we showed above that the functions
$h_{jk}$ are linearly independent.  The same reasoning we used there
yields that such a relation can only be true if
$$
\sum_{j=1}^{n-1} c_{ij}S(a_j,a_k)=
\cases
1, &\text{if }i=k \\
0, &\text{if }i\ne k.
\endcases
$$
This shows that the $(n-1)\times(n-1)$ matrix $\left[S(a_j,a_k)\right]$
is invertible and $[c_{ij}]$ is its inverse.

We remark here that formula (3.2) has an interesting application.
It is quite easy to prove that $S(z,a)$ is in $C^\infty(\Obar)$ as a function
of $z$ for each fixed $a\in\O$ (see \cite{2, page~22}).  The more difficult
result that $S(z,w)$ is in $C^\infty(\Obar\times\Obar-\{(z,z):z\in b\O\})$
follows directly from the the smoothness of $S(z,a)$ for fixed $a$ and
formula (3.2).

\subhead 4. Complexity of the Szeg\H o kernel \endsubhead
Formula (3.2) reveals that the Szeg\H o kernel associated to an
$n$-connected domain is composed of the $n+1$ functions, $S(z,a)$,
$S(z,a_1)$, $S(z,a_2),\dots,S(z,a_{n-1}),$ and $L(z,a)$ (because
$f(z)=S(z,a)/L(z,a)$).  If one knows the boundary values of these $n+1$
functions, then the Szeg\H o kernel may be evaluated at any pair of points
$(z,w)$ in $\O\times\O$ by applying the Cauchy integral formula twice,
once to evaluate the functions on the right hand side of (3.2) at $z$, and
once to evaluate the functions at $w$.  In this section, we show how much
effort is required to numerically compute the boundary values of the $n+1$
functions comprising $S(z,w)$.

Kerzman and Stein \cite{10} discovered an effective method for
computing the Szeg\H o kernel (see also \cite{{2,3,7,11,14}}).  They proved
that the function $\Cal S_a(z)=S(z,a)$ is the solution to a Fredholm
integral equation of the second kind given by
$$\Cal S_a(z)-\int_{w\in b\O}A(z,w)\Cal S_a(w)\ ds=\Cal C_a(z),$$
where $A(z,w)$ is the Kerzman-Stein kernel and $\Cal C_a(z)$ is the Cauchy
kernel.  To be precise,
$$A(z,w)=\frac{1}{2\pi i}\left(
\frac{T(w)}{w-z}-\frac{\overline{T(z)}}{\bar w-\bar z}\right)$$
if $z,w\in b\O$, $z\ne w$, and $A(z,w)=0$ if $z=w$, and
$$\Cal C_a(z)=\frac{1}{2\pi i}\frac{\overline{T(z)}}{\bar a-\bar z}.$$
The Kerzman-Stein kernel is skew-hermitian and, in spite of the
apparent singularity at $z=w$ in the formula above, it is in
$C^\infty(b\O\times b\O)$.  (Kerzman and Stein discovered that the apparent
singularities in the formula for $A(z,w)$ exactly cancel.) \ The
Cauchy kernel is in $C^\infty(b\O)$.  It follows from standard
theory that this integral equation has a unique $C^\infty$ smooth
solution.  (See Kerzman and Tummer \cite{14} and \cite{{3,7}} for
descriptions of convenient ways to write and to solve this integral
equation.)

The Kerzman-Stein equation produces the boundary values of $S(z,a)$.
The boundary values of the Garabedian kernel $L(z,a)$ can be computed
via identity (2.1), and the boundary values of the Ahlfors map $f_a(z)$
can now be gotten from (2.2).  The remaining functions in expansion (3.2) can
be computed via the Kerzman-Stein integral equation once the zeroes
$a_1,\dots,a_{n-1}$ have been located.  Since $S(z,a)$ does not vanish
on $b\O$, we may use the residue theorem to compute the symmetric sums
$$\sum_{j=1}^{n-1}a_j^k=\int_{z\in b\O}
\frac{z^k(\dee/\dee z)S(z,a)}{S(z,a)} \ dz$$
for $k=1,\dots,n-1$.  Newton's identities can now be used to compute
the elementary symmetric functions of $a_1,\dots,a_{n-1}$, and hence, the
coefficients of the polynomial $\prod_{j=1}^{n-1}(\zeta-a_j)$ are
determined.  We have therefore shown that the problem of locating the
zeroes of $S(z,a)$ is equivalent to computing $n-1$ line integrals
and finding the roots of a polynomial of degree $n-1$.

\subhead 5. The Bergman kernel \endsubhead
In this section, we shall prove that the Bergman kernel of an
$n$-connected domain in the plane with $C^\infty$ smooth boundary
is composed of the same basic functions that comprise the Szeg\H o kernel.
We shall also prove that the Bergman kernel can be computed at every
pair of points by solving $n$ one dimensional $C^\infty$ Fredholm integral
equations of the second kind, and by solving a linear system.  At no point
is it necessary to evaluate a double integral.

The Bergman kernel $K(z,w)$ is related to the Szeg\H o kernel via the
identity
$$K(z,w)=4\pi S(z,w)^2+\sum_{i,j=1}^{n-1}
A_{ij}F_i'(z)\overline{F_j'(w)},$$
where the functions $F_i'(z)$ are classical functions of potential theory
described as follows.  The harmonic function $\omega_j$ which solves the
Dirichlet problem on $\O$ with boundary values one on the boundary
curve $\gamma_j$ and zero on $\gamma_k$ if $k\ne j$ has a multivalued
harmonic conjugate.  The function $F_j'(z)$ is a globally defined
single valued holomorphic function on $\O$ which is locally defined as
the derivative of $\omega_j+iv$ where $v$ is a local harmonic
conjugate for $\omega_j$.  The Cauchy-Riemann equations reveal that
$F_j'(z)=2(\dee\omega_j/\dee z)$.

Let $\Cal F'$ denote the vector space of functions given by the complex
linear span of the set of functions $\{F_j'(z):j=1,\dots,n-1\}$.  It is
a classical fact that $\Cal F'$ is $n-1$ dimensional.  Notice that
$S(z,a_i)L(z,a)$ is in $A^\infty(\O)$ because the pole
of $L(z,a)$ at $z=a$ is cancelled by the zero of $S(z,a_i)$ at $z=a$.
A theorem due to Schiffer (see \cite{{12,2,4}}) states that the $n-1$
functions $S(z,a_i)L(z,a)$, $i=1,\dots,n-1$ form a basis for $\Cal F'$.  We
may now write
$$K(z,w)=4\pi S(z,w)^2+\sum_{i,j=1}^{n-1}
\lambda_{ij}S(z,a_i)L(z,a)\,\overline{S(w,a_j)L(w,a)},\tag5.1$$
which, together with (3.2) allows us to write down a formula which sheds
light on the degree of complexity of the Bergman kernel.
\proclaim{Theorem 5.1}
The Bergman kernel is composed of the same basic functions that make
up the Szeg\H o kernel, as evidenced by the following formula.
$$\gather
K(z,w)=\frac{1}{(1-f(z)\overline{f(w)})^2}\left(
\sum\Sb
0\le i\le j\le n-1 \\
0\le k\le m\le n-1
\endSb
C_{ijkm}S(z,a_i)S(z,a_j)\,\overline{S(w,a_k)S(w,a_m)}\right) \\
+\sum_{i,j=1}^{n-1}\lambda_{ij}S(z,a_i)L(z,a)\,\overline{S(w,a_j)L(w,a)}.
\endgather$$
\endproclaim

We shall now discuss the amount of computational effort required to compute
the Bergman kernel.  To streamline the argument, it will be convenient
to use the fact that the linear span of
$\{S(z,a_i)L(z,a): i=1,\dots,n-1\}$ is the same as the linear span
of $\{L(z,a_i)S(z,a): i=1,\dots,n-1\}$ (see \cite{2, page~80}).  Hence,
formula (5.1) can also be written in the form
$$K(z,w)=4\pi S(z,w)^2+\sum_{i,j=1}^{n-1}
\lambda_{ij}L(z,a_i)S(z,a)\,\overline{L(w,a_j)S(w,a)},\tag5.2$$
(where, here, the coefficients $\lambda_{ij}$ represent different
constants than they do in (5.1)).
The difficulty of computing the functions appearing in (5.2) has been
discussed.  We now describe a method for computing the coefficients
$\lambda_{ij}$.  We shall write $K_w(z)$ in place of $K(z,w)$ and $S_w(z)$
in place of $S(z,w)$ to emphasize that we are thinking of $w$ as being
fixed and we are viewing these kernels as functions of $z$.  Let us also
write $\Cal L_i(z)=L(z,a_i)S(z,a)$.  Thus, formula (5.2) may be rewritten as
$$K_w-4\pi S_w^2=\sum_{i,j=1}^{n-1}
\lambda_{ij}\overline{\Cal L_j(w)}\,\Cal L_i.\tag5.3$$

The complement of $\Obar$ in $\C$ is the union of domains $D_j$,
$j=1,\dots,n$, where the boundary of $D_j$ is described by the
boundary curve $\gamma_j$.  Recall that $\gamma_n$ denotes the outer
boundary curve of $\O$.  For $j=1,\dots,n-1$, pick a point $b_j$ in $D_j$.
We now consider the effect of integrating (5.3) against the function
$1/(z-b_k)$.  Notice that
$$\langle (z-b_k)^{-1},K_w\rangle_{\O}=\frac{1}{w-b_k}$$
because the Bergman kernel reproduces holomorphic functions.
Since $1/(z-b_k)=(\dee/\dee z)\ln|z-b_k|^2$, we may use
the complex Green's identity to compute
$$\gather
\langle (z-b_k)^{-1},S_w^2\rangle_{\O}=
\iint_{z\in\O}
(\dee/\dee z)\ln|z-b_k|^2
\ \overline{S_w(z)^2}\ (\frac{i}{2}dz\wedge d\bar z)= \\
i\int_{z\in b\O}\ln|z-b_k|\ \overline{S_w(z)^2}\ d\bar z.
\endgather$$
Define numbers $A_{ik}= \langle (z-b_k)^{-1},\Cal L_i\rangle_{\O}$.
We may use the complex Green's identity again to obtain
$$A_{ik}=
i\int_{z\in b\O}\ln|z-b_k|\ \overline{\Cal L_i(z)}\ d\bar z.$$
We now collect the integrals above as dictated by (5.3), and
we set $w=a_m$, $m=1,\dots,n-1$ to obtain the system,
$$\frac{1}{a_m-b_k}-
4\pi i\int_{z\in b\O}\ln|z-b_k|\ S(a_m,z)^2\ d\bar z=\sum_{i,j=1}^{n-1}
\lambda_{ij}A_{ik}\overline{\Cal L_j(a_m)}.$$
To show that this system determines the numbers $\lambda_{ij}$, we
need only check that the matrices given by $\Bbb A=[A_{ik}]$ and
$\Bbb L=[\Cal L_j(a_m)]$ are invertible.  That $\Bbb L$ is invertible
is obvious because
$$L(w,a_j)S(w,a)=\cases
0, &\text{if $w=a_m$, $m\ne j$} \\
\frac{1}{2\pi}\frac{\dee}{\dee z}S(a_j,a), &\text{if $w=a_j$},
\endcases$$
and $(\dee/\dee z)S(a_j,a)\ne0$ because $a$ has been chosen so that the
zeroes of $S(z,a)$ are simple zeroes.
To show that $\Bbb A$ is invertible, we shall need to use an argument
from \cite{4}.  If $G=\sum_{k=1}^{n-1}c_kF_k'$, then $G=2(\dee/\dee z)\omega$
where $\omega=\left(\sum_{k=1}^{n-1}c_k\omega_k\right)$.  It is
proved in \cite{4, page~12} that the constants $c_k$ are given by the integral
$$c_k=-\frac{1}{2\pi i}\int_{z\in b\O}\ln|z-b_k|G(z)dz,$$
where $b_k$ is the fixed point chosen from $D_k$.  Notice that $c_k$
is the value of $\omega$ on $\gamma_k$.
Suppose $\Bbb A$ is not invertible.  Then there would exist constants
$\sigma_i$, not all zero, such that
$$\sum_{i=1}^{n-1}A_{ik}\bar\sigma_i=0$$
for each $k$, and the complex conjugate of this equality yields that
$$\int_{z\in b\O}\ln|z-b_k|\left(
\sum_{i=1}^{n-1}\sigma_i\Cal L_i(z)\right)\ dz=0.\tag5.4$$
Let $G=\sum_{i=1}^{n-1}\sigma_i\Cal L_i$.  Since $G$
is in the linear span of $\{F_j'\}_{j=1}^{n-1}$, condition (5.4)
and the fact from \cite{4} imply that $G=2(\dee/\dee z)\omega$ where
$\omega$ is a harmonic function on $\O$ that vanishes on each boundary
curve of $\O$, i.e., that $G\equiv0$.  Now each $\sigma_k$
must be zero because the functions $\Cal L_i$ are linearly independent.
This contradiction yields that the matrix $\Bbb A$ must be non-singular and
the proof is finished.

\subhead 6. Characterization of domains with rational kernel functions
\endsubhead
In the previous sections, we have shown that the kernel functions are
not as complex as one might expect them to be.  In this section,
we shall prove theorems that say roughly that the only domain whose
Bergman or Szeg\H o kernels are so simple as to be rational functions
is the disc.

A function $R(z,w)$ of two complex variables is called rational if there
are relatively prime polynomials $P(z,w)$ and $Q(z,w)$ such that
$R(z,w)=P(z,w)/Q(z,w)$.   It is not hard to prove that a function $H(z,w)$,
which is holomorphic in $z$ and $w$ on a product domain $\O_1\times\O_2$ is
rational if and only if, for each fixed $b\in\O_2$, the function $H(z,b)$ is
rational in $z$, and for each fixed $a\in\O_1$, the function $H(a,w)$ is
rational in $w$ (see Bochner and Martin \cite{6, page~201}).  We shall say
that the Bergman kernel function $K(z,w)$ associated to a domain $\O$ is
rational if it can be written as
$R(z,\bar w)$ where $R$ is a holomorphic rational function of two
variables.  Because the Bergman kernel is hermitian, the facts above imply
that $K(z,w)$ is rational if and
only if, for each point $a\in\O$, the function $K(z,a)$ is a rational
function of $z$.  In fact, $K(z,w)$ is rational if and only if there
exists a small disc $D_\epsilon(w_0)\subset\O$ such that $K(z,a)$ is
a rational function of $z$ for each $a\in D_\epsilon(w_0)$.  Similar
statements hold for the other kernel functions.

\proclaim{Theorem 6.1}
Suppose $\O$ is a bounded $n$-connected domain, $n>1$, with $C^\infty$ smooth
boundary.  Neither the Bergman kernel nor the Szeg\H o kernel associated
to $\O$ can be rational functions.
\endproclaim

The assumption in Theorem~6.1 that the boundary of $\O$ is $C^\infty$
smooth can be relaxed.  For example, the conclusion about the Szeg\H o
kernel holds if the boundary is only assumed to be $C^2$ smooth.  The
conclusion about the Bergman kernel holds if the domain is only assumed
to be finitely connected and such that no boundary component is a point.
We shall explain how to relax the smoothness assumptions later in this section.

Before we proceed to prove Theorem~6.1, let us consider the case of a
{\it one\,}-connected domain $\O\ne\C$.  If $f_a$ is a Riemann mapping
$f_a:\O\to D_1(0)$ such that $f_a(a)=0$ and $f_a'(a)>0$, the Bergman kernel
for $\O$ can be expressed via
$$K(z,w)=\frac{f_a'(z)\,\overline{f_a'(w)}}
{\pi(1-f_a(z)\,\overline{f_a(w)}\ )^2}.$$
If we set $w=a$ in this formula, we obtain the identity
$$K(z,a)=C\,f_a'(z),$$
where $C=f_a'(a)/\pi$ is a positive constant.
If we differentiate the formula with respect to $\bar w$ and then set $w=a$,
we obtain
$$\frac{\dee}{\dee\bar w}K(z,a)=f_a'(z)(C_1+C_2f_a(z)),$$
where $C_1$ and $C_2$ are constants, and $C_2\ne0$.  (In fact,
$C_2=2f_a'(a)^2/\pi$.) \ It can easily be deduced from these formulas
that the Bergman kernel is rational if and only if the Riemann map is
rational.

To study the Szeg\H o kernel, assume that $\O$ is a bounded simply
connected domain with $C^2$ smooth boundary and let $f_a$ denote a
Riemann map as above.  The Szeg\H o kernel is given by
$$S(z,w)=\frac{\sqrt{f_a'(z)}\ \overline{\sqrt{f_a'(w)}}}
{2\pi(1-f_a(z)\,\overline{f_a(w)}\ )}.$$
Set $w=a$ in this formula to obtain
$$S(z,a)=c\sqrt{f_a'(z)},$$
where $c=\sqrt{f_a'(a)}/(2\pi)$ is a positive constant.  Now differentiate the
formula with respect to $\bar w$ and then set $w=a$ to obtain
$$\frac{\dee}{\dee\bar w}S(z,a)=\sqrt{f_a'(z)}(c_1+c_2f_a(z)),$$
where $c_1$ and $c_2$ are constants, and $c_2\ne0$.  (In fact,
$c_2=f_a'(a)^{3/2}/(2\pi)$.) \ These formulas reveal that Szeg\H o
kernel is rational if and only if the Riemann map and the square root
of its derivative are rational.  Let us summarize these results in the
following theorem

\proclaim{Theorem 6.2}
Suppose $\O\ne\C$ is a simply connected domain.  The Bergman kernel
associated to $\O$ is rational if and only if there is a rational
biholomorphic mapping $f(z)$ mapping $\O$ one-to-one onto the unit disc.
If $\O$ is further assumed to be bounded and have $C^2$ smooth
boundary, then the Szeg\H o kernel associated to $\O$ is rational if and only
if there is a rational biholomorphic mapping $f(z)$ mapping $\O$ one-to-one
onto the unit disc such that $f'(z)$ is the square of a rational
function.
\endproclaim

\demo{Proof of Theorem 6.1}
We shall use the notations that we set up previously to describe our
$n$-connected domain $\O$.  Hence, $\gamma_n$ denotes the outer
boundary of $\O$.  Since we are assuming that $n>1$, we may let
$\gamma_1$ denote one of the inner boundary curves of $\O$, and we
let $D_1$ denote the bounded region enclosed by $\gamma_1$.

We first assume that the Szeg\H o kernel associated to $\O$ is rational.
Formula (3.2) shows that it then follows that the Ahlfors mapping
$f_a(z)$ is a rational function of $z$ for each point $a$ in $\O$
minus the finite set where the zeroes of $S(z,a)$ might not all be
simple zeroes.  (It was proved earlier by M.~Jeong \cite{{8,9}}, using
other techniques, that the Ahlfors maps are all rational.) \ We may
now use formula (2.2) to deduce that the Garabedian kernel $L(z,a)$
is a rational function of $z$ for each $a$ in an open subset of $\O$, and
hence that $L(z,a)$ is a rational function of $(z,a)$.  It is clear that the
Ahlfors maps $f_a$ can have no poles on $b\O$.  Since, the boundary
of $\O$ is assumed to be smooth, the Hopf lemma implies that $f_a'(z)\ne0$
for $z\in b\O$.  Since the boundary curves of $\O$ are described by
the equation $|f_a(z)|$=1, it follows that the boundary curves of $\O$
are all {\it real analytic curves}.  From this it follows that
$S(z,w)$ extends holomorphically in $z$ and antiholomorphically in $w$ to an
open set in $\C\times\C$ containing $\Obar\times\Obar-\{(z,z):z\in b\O\}$,
and that $L(z,w)$ extends holomorphically in $z$ and $w$ to an
open set in $\C\times\C$ containing $\Obar\times\Obar-\{(z,z):z\in\Obar\}$.
This may seem like a silly thing to say in light of the fact that
$S(z,w)$ and $L(z,w)$ are rational, however it implies that singularities
must stay away from $b\O\times b\O-\{(z,z):z\in b\O\}$.

We shall be concerned with the number of zeroes and poles of $S(z,a)$
and $L(z,a)$ as functions of $z$ which lie in $D_1$, and we shall consider how
these numbers vary as $a$ moves from a point on the outer boundary of $\O$
to a point on $\gamma_1$.  First, however, we shall need to review some
properties of the zeroes of the Szeg\H o kernel proved in \cite{4} (see also
\cite{2}).  We mentioned earlier that if $a\in\O$, then $S(z,a)\ne0$ and
$L(z,a)\ne0$ for all $z\in b\O$.  We also mentioned that neither $S(z,a)$ nor
$L(z,a)$ can have poles on $b\O$.  We shall use these facts to see that zeroes
and poles of $S(z,a)$ and $L(z,a)$ which lie in $D_1$ cannot exit $D_1$ through
$\gamma_1$ as $a$ varies in $\O$.  We also mentioned earlier that $S(z,a)$ has
$n-1$ zeroes in $\O$ as a function of $z$, and $L(z,a)\ne0$ for all
$z\in\Obar-\{a\}$.  It
is proved in \cite{4} that, if $a\in\O$ is allowed to tend to a point $A_k$
in a boundary curve $\gamma_k$, then the $n-1$ zeroes of $S(z,a)$
separate into simple zeroes which migrate to distinct points on the
boundary in such a way that there is a point on each boundary curve
$\gamma_j$, $j\ne k$ to which exactly one of the zeroes tends.  To be
precise, there exist points $\{A_j: 1\le j\le n,\ j\ne k\}$
with $A_j\in\gamma_j$ such that the $n-1$ zeroes of $S(z,a)$ can be listed as
$a_j$, $1\le j\le n$, $j\ne k$, where $a_j$ tends to $A_j$ for each $j\ne k$
as $a$ tends to $A_k$.

Since $S(z,w)$ is rational, there exist relatively prime polynomials
$P(z,w)$ and $Q(z,w)$ such that $S(z,w)=P(z,\bar w)/Q(z,\bar w)$.  There
are at most finitely many points $w_0\in\C$ for which the equations
$P(z,\bar w_0)=0$ and $Q(z,\bar w_0)=0$ have a common root (see Ahlfors
\cite{1, page~300}).  Let $B_S$ denote the (possibly empty) set of such
points $w_0$.  Similarly, there are relatively prime polynomials $p(z,w)$
and $q(z,w)$ such that $L(z,w)=p(z,w)/q(z,w)$, and there is a finite set
$B_L$ of points $w_0$ where the equations $p(z,w_0)=0$ and $q(z,w_0)=0$ have
a common root.  Let $B=B_S\cup B_L$.

Let $S_a(z)=S(z,a)$.  It is a simple exercise using the argument principle
that the zeroes and poles of $S_a$ are continuous
functions of $a$ when $a\not\in B$ in the following sense.  Suppose $z_0$
is a zero of multiplicity $m$ of $S(z,a_0)$ where $a_0\not\in B$.
Given $\epsilon>0$ such that $z_0$ is the only zero of $S(z,a_0)$ in
$\overline{D_\epsilon(z_0)}$, there is a $\delta>0$ such that
$S(z,a)$ has precisely $m$ zeroes in $\overline{D_\epsilon(z_0)}$
as a function of $z$ (counting multiplicities) when $a\in D_\delta(a_0)$.
A similar statement holds for poles of $S(z,a)$, and for zeroes and
poles of $L(z,a)$.

We have stated all the necessary facts to be able to assert that there
exist non-negative integers $Z_S$, $Z_L$, $P_S$,  and $P_L$ such that,
for any point $a\in\O-B$,
$Z_S$ is equal to the number of zeroes of $S(z,a)$ in $\overline{D_1}$,
$Z_L$ is equal to the number of zeroes of $L(z,a)$ in $\overline{D_1}$,
$P_S$ is equal to the number of poles of $S(z,a)$ in $\overline{D_1}$, and
$P_L$ is equal to the number of poles of $L(z,a)$ in $\overline{D_1}$.

Let $\sigma$ denote a curve in $\Obar-B$ which starts at a point $A_n$
on the outer boundary $\gamma_n$ of $\O$, travels through $\O$, and
terminates at a point $A_1$ in $\gamma_1$.  We shall be able to deduce
relationships between the four integers, $Z_S$, $Z_L$, $P_S$,  and $P_L$,
by letting $a$ tend to the two endpoints of $\sigma$.  The
relationships shall turn out to be contradictory.  To find
relationships between these numbers, we shall need to use an argument
from \cite{4}.  Since the boundary curves of $\O$ are real analytic
curves, there exists an {\it antiholomorphic reflection function\/}
$R(z)$ with the properties that $R(z)$ is defined and is antiholomorphic
on a neighborhood $\Cal O$ of $b\O$, $R(z_0)=z_0$ when $z_0\in b\O$,
$R'(z)$ is non-vanishing on $\Cal O$, and $R(z)$ maps $\Cal O\cap\O$
one-to-one onto $\Cal O-\Obar$.

Let $w_k$ be a sequence of points in $\O$ that tend to $A_n$ along $\sigma$,
and let $a$ be a fixed point in $\O-B$.  By (2.1), we have
$-i\,L(z,a)T(z)=S(a,z)$ and $-i\,L(z,w_k)T(z)=S(w_k,z)$ for $z\in b\O$.
Divide the second of these identities by the first and use the fact
that $R(z)=z$ on $b\O$ to obtain
$$\frac{S(w_k,z)}{S(a,z)}=\frac{L(R(z),w_k)}{L(R(z),a)}
\qquad\text{for $z\in b\O$}.\tag6.1$$
The function on the left hand side of (6.1) is antiholomorphic in $z$
on a neighborhood of $b\O$; so is the function on the right hand side.
Since these functions agree on $b\O$, they must be equal on a
neighborhood of $b\O$.  In fact, because $S$ and $L$ are rational,
these two functions are equal as meromorphic functions on the
neighborhood $\Cal O$ of $b\O$ on which $R(z)$ is defined.  We may
assume that $\Cal O$ is small enough that $S(z,a)$ and $L(z,a)$ have
no poles or zeroes in $\Cal O$.  Formula (6.1) now allows us to read
off the following facts (keep in mind that $w_k$ is close to
$A_n\in\gamma_n$).  
If $S(z,w_k)$ has a zero $z_0\in\O$ near $\gamma_1$, then $L(R(z_0),w_k)=0$,
i.e., $L(z,w_k)$ has a zero at the reflected point $R(z_0)\in D_1$ near
$b\O$.  Neither $S(z,w_k)$ nor $L(z,w_k)$ can have a  pole $z_0\in D_1$
near $\gamma_1$, because neither $L(z,w_k)$ nor $S(z,w_k)$ has a pole at the
reflected point $R(z_0)\in\O$.

Finally, notice that (2.1) yields that
$$-i\,L(A_n,z)T(A_n)=S(z,A_n)\qquad\text{for $z\in\O$,}$$
and consequently $-i\,L(A_n,z)T(A_n)=S(z,A_n)$ for $z\in\overline{D_j}$.
Hence, the functions $L(A_n,z)$ and $S(z,A_n)$ have the same number of
zeroes and poles in $\overline{D_j}$.  It is proved in \cite{4} that
$S(z,A_n)$ has a single simple zero on $\gamma_1$ and this zero is
approached by single simple zeroes of $S(z,w_k)$ as $k\to\infty$.  No
other zeroes of $S(z,w_k)$ can migrate near $\gamma_1$.  Our remarks above
yield that $L(z,w_k)$ has a simple zero at the reflection of the zero of
$S(z,w_k)$ near $\gamma_1$.  By letting $k\to\infty$, we obtain the relations
$$\align
Z_L &= Z_S+1 \\
P_L &= P_S.
\endalign$$

We now take a sequence of points $w_k$ in $\O$ that tend to $A_1$ along
$\sigma$.  Formula (6.1) remains valid and, if we reason as above, we deduce
that, since $S(z,w_k)$ has no zeroes $z_0$ with $z_0$ near $\gamma_1$,
$L(z,w_k)$ has no zeroes in $D_j$ near $\gamma_1$.  However, since
$L(z,w_k)$ has a simple pole at $z=w_k$, it follows that
$S(z,w_k)$ has a simple pole at the reflected point $R(w_k)$ in $D_1$.
We now let $k\to\infty$ and use the facts that $S(z,A_1)$ and $L(z,A_1)$
have the same zeroes and poles in $\overline{D_1}$ and that one of those
poles is a simple pole at $A_1$ to obtain
$$\align
Z_L &= Z_S \\
P_L+1 &= P_S.
\endalign$$
These relationships contradict the ones we obtained by letting $w_k$
tend to $A_n$, and we conclude that $\O$ cannot be multiply
connected.

We now turn to the study of the Bergman kernel.  Assume that $K(z,w)$
is rational.  We shall use an argument similar to the one above for
the Szeg\H o kernel, however, many of the underlying facts are
different.  Before we can begin, we must review some facts about the
Bergman kernel (see \cite{2} for proofs of these facts).

We first must prove that if the Bergman kernel associated to a
bounded domain is rational, then any proper holomorphic mapping of the
domain onto the unit disc must be rational.  Suppose $f:\O\to D_1(0)$
is a proper holomorphic map.  Such a map must be in $A^\infty(\O)$ and
there is a positive integer $m$ such that $f$ is an $m$-to-one mapping
of $\O$ onto $D_1(0)$ (see \cite{2, page~62--70}).  The branch locus
$\Cal B=\{z\in\O: f'(z)=0\}$ is a finite set, and for each point
$w_0$ in $D_1(0)-f(\Cal B)$, there are exactly $m$ distinct points in
$f^{-1}(w_0)$.  Near such a point $w_0$, there is an $\epsilon>0$ such
that it is possible to define $m$ holomorphic maps $F_1(w),\dots,F_m(w)$
on $D_\epsilon(w_0)$ which map into $\O-\Cal B$ such that $f(F_k(w))=w$.
These local inverses appear in the following transformation formula for the
Bergman kernels under a proper holomorphic mapping.  Let
$K_1(z,w)=\pi^{-1}(1-z\bar w)^{-2}$ denote the Bergman kernel of the unit
disc (and recall that $K(z,w)$ denotes the Bergman kernel for $\O$).  It
is proved in \cite{2, page~68} that the kernels transform via
$$f'(z)K_1(f(z),w)=\sum_{k=1}^m K(z,F_k(w))\overline{F_k'(w)}.$$
Although the functions $F_k$ are only locally defined on $D_1(0)-f(\Cal B)$,
the function on the right hand side of the transformation formula, being
symmetric in the $F_k$, is globally well defined.  In fact, the
function is holomorphic in $z$ and antiholomorphic in $w$ for
$(z,w)\in\O\times(D_1(0)-f(\Cal B))$.  (The set $f(\Cal B)$
can be seen to be a removable singularity set, but we shall not need
to know this.) If the origin is in $D_1(0)-f(\Cal B)$, we replace $f$ by
its composition with a M\"obius transformation so that $0\not\in
D_1(0)-f(\Cal B)$.  We now set $w=0$ in the transformation formula for the
Bergman kernels to obtain
$$f'(z)=\pi\sum_{k=1}^m K(z,F_k(0))\overline{F_k'(0)}.$$
This shows that $f'(z)$ is a rational function.  Now differentiate the
transformation formula with respect to $\bar w$ and then set $w=0$
to obtain
$$2f'(z)f(z)=\pi\sum_{k=1}^m
\frac{\dee}{\dee\bar w}K(z,F_k(0))\overline{F_k'(0)^2}
+\pi\sum_{k=1}^m K(z,F_k(0))\overline{F_k''(0)}.$$
We may now deduce that $f'(z)f(z)$ is rational, and so it follows that
$f(z)$ is rational.

Since the Ahlfors mappings $f_a(z)$ are proper mappings of $\O$ onto the unit
disc, they are rational functions of $z$.  As above, this implies that the
boundary curves of $\O$ are all {\it real analytic curves}, and from this
it follows that $K(z,w)$ extends holomorphically in $z$ and
antiholomorphically in $w$ to an open set in $\C\times\C$ containing
$\Obar\times\Obar-\{(z,z):z\in b\O\}$.

The Bergman kernel is related to the classical Green's function via
$$K(z,w)=-\frac{2}{\pi}\frac{\dee^2 G(z,w)}{\dee z\dee\bar w}.$$
Define another function $\Lambda(z,w)$ on $\O$ via
$$\Lambda(z,w)=-\frac{2}{\pi}\frac{\dee^2 G(z,w)}{\dee z\dee w}.$$
(This function is sometimes written $L(z,w)$ in the
literature; we have chosen the symbol $\Lambda$ here to avoid
confusion with our notation for the Garabedian kernel above.)
It follows from known properties of the Green's function that
$\Lambda(z,w)$ extends holomorphically in $z$ and $w$ to an
open set in $\C\times\C$ containing $\Obar\times\Obar-\{(z,z):z\in\Obar\}$,
and that, if $a\in\Obar$, then $\Lambda(z,a)$ has a double pole at $z=a$
as a function of $z$.

We shall need to use the following real variable theorem.
Suppose that $R(x,y)$ is a real analytic function of $(x,y)$ on a
product domain $U_1\times U_2 \subset\R^n\times\R^m$ such that
$R(x,y_0)$ is a rational function of $x$ on $U_1$ for each $y_0\in U_2$, and
$R(x_0,y)$ is a rational function of $y$ on $U_2$ for each $x_0\in U_1$.
It then follows that $R(x,y)$ is a rational function of $(x,y)$ (the
proof in Bochner and Martin \cite{6, page~201} works in the real case,
too).

We now consider the function $|f_a(z)|^2=|S(z,a)|^2/|L(z,a)|^2$.  We
know that for each fixed $a\in\O$, the function  $f_a(z)$ is a rational
function of $z$, and hence $|f_a(z)|^2$ is a rational function of $(x,y)$
where $z=x+iy$.  Since $|f_a(z)|=|f_z(a)|$, the real variable theorem
mentioned above implies that $|f_a(z)|^2$ is a rational function
of the variables $\text{Re }z$, $\text{Im }z$, $\text{Re }a$, and
$\text{Im }a$.

Assume that $a\in\O$ is such that the zeroes of $S(z,a)$ are all simple zeroes.
Because Ahlfors maps are proper holomorphic maps, it is easy to verify that
$$\frac{1}{2}\ln |f_a(z)|^2=G(z,a)+\sum_{i=1}^{n-1}G(z,a_i)\tag6.2$$
where the points $a_i$, $i=1,\dots,n-1$ are the zeroes of $S(z,a)$ (which,
together with $a$, are the zeroes of $f_a$).  We now consider the
way in which the zeroes $a_i$ depend on $a$, and we write
$a_i(a)$ in order to regard $a_i$ as a function of $a$.  Let $A_0$ be
a fixed point in $\O$ such that the zeroes of $S(z,A_0)$ are simple.
Since the points $a_i(A_0)$ are distinct, we may choose an
$\epsilon>0$ such that $\overline{D_\epsilon(a_i(A_0))}\subset\O$ for each
$i$ and $\overline{D_\epsilon(a_i(A_0))}\cap\overline{D_\epsilon(a_j(A_0))}
=\emptyset$ if $i\ne j$.  Thus, $a_i(A_0)$ is the only zero of $S(z,A_0)$ in
$\overline{D_\epsilon(a_i(A_0))}$.  The dependence
of the zeroes of $S(z,a)$ on $a$ can be described by the formula,
$$a_i(a)=\frac{1}{2\pi i}\int_{|z-a_i(A_0)|=\epsilon}
z\,\frac{\frac{\dee}{\dee z}S(z,a)}{S(z,a)}\ dz,$$
which is valid when $a$ is close to $A_0$.
Because $S(z,a)$ is antiholomorphic in $a$, this formula shows that
$a_i(a)$ is an antiholomorphic function of $a$ near $A_0$.  We now
differentiate (6.2) with respect to $z$ to obtain
$$\frac{f_a'(z)}{2f_a(z)}=\frac{\dee}{\dee z} G(z,a)+
\sum_{i=1}^{n-1} \frac{\dee}{\dee z} G(z,a_i).$$
Next, we differentiate with respect to $a$ and use the complex chain
rule to obtain
$$\frac{\dee}{\dee a}\left(\frac{f_a'(z)}{2f_a(z)}\right)
=\frac{\dee^2 G(z,a)}{\dee z\dee a}+
\sum_{i=1}^{n-1} \frac{\dee^2 G(z,a_i)}{\dee z\dee\bar a_i}
\frac{\dee\bar a_i}{\dee a}.\tag6.3$$
We now claim that the function on the left hand side of (6.3) is
is a rational function $R(z,a)$ of $z$ and $a$.  Indeed,
because $|f_a(z)|^2$ is rational in the real and imaginary parts of
$z$ and $a$, it follows that $R(z,a)$ is rational in the real and
imaginary parts of $z$ and $a$.  It is clear that $R(z,a)$ is
holomorphic in $z$.  Since $|f_a(z)|=|f_z(a)|$, it follows that
$R(z,a)=R(a,z)$, and so $R(z,a)$ is holomorphic in $a$, too.
Consequently, $R(z,a)$ is a rational function of $z$ and $a$.
The function on the right hand side of (6.3) can be rewritten to yield
$$R(z,a)=-\frac{\pi}{2}\Lambda(z,a)-\frac{\pi}{2}
\sum_{i=1}^{n-1}K(z,a_i)\frac{\dee\bar a_i}{\dee a}.$$
This last formula shows that, for each fixed $a$ in an open subset of
$\O$, the function $\Lambda(z,a)$ is a rational function of $z$.  Since
$\Lambda(z,a)=\Lambda(a,z)$, we conclude that $\Lambda(z,a)$ is a
rational function of $(z,a)$.

The Bergman kernel is related to $\Lambda$ via the identity
$$\Lambda(w,z)T(z)=-K(w,z)\overline{T(z)}\qquad\text{for $w\in\O$ and
$z\in b\O$}\tag6.4$$
(see \cite{2, page~135}).  Define the set $B$ to be the finite set of points
$a$ at which the numerators and denominators of $K(z,a)$ and $\Lambda(z,a)$
have common zeroes as functions of $z$.  Let $\sigma$ denote a curve in
$\Obar-B$ which starts at a point $A_n$ on the outer boundary $\gamma_n$ of
$\O$, travels through
$\O$, and terminates at a point $A_1$ on an inner boundary curve $\gamma_1$.
Since $K(z,a)$ and $\Lambda(z,a)$ cannot have poles on the boundary
as functions of $z$ when $a\in\O$, the number $P_K$ of poles of
$K(z,a)$ as a function of $z$ which lie in $D_1$ is constant as $a$ moves
along $\sigma$ away from the endpoints of the curve.  Also,
the number $P_\Lambda$ of poles of $\Lambda(z,a)$ in $D_1$ is constant
as $a$ moves along the curve.  We shall deduce relationships between $P_K$
and $P_\Lambda$ by letting $a$ tend to the two endpoints of $\sigma$.  The
relationships shall turn out to be contradictory.

Let $w_k$ be a sequence of points in $\O$ that tend to $A_n$ along $\sigma$.
Since $K(z,a)$ and $\Lambda(z,a)$ cannot vanish for $z\in b\O$ when
$a$ is close to the boundary (see \cite{13} or \cite{2, page~132}), we
may choose a point $a$ in $\O-B$ so that these functions are
non-vanishing in $z$ near $b\O$.  By (6.4), we have
$\Lambda(a,z)T(z)=-K(a,z)\overline{T(z)}$ and
$\Lambda(w_k,z)T(z)=-K(w_k,z)\overline{T(z)}$ for $z\in b\O$.
Divide the second of these identities by the first and use the fact
that $R(z)=z$ on $b\O$ to obtain
$$\frac{\Lambda(w_k,z)}{\Lambda(a,z)}=\frac{K(w_k,R(z))}{K(a,R(z))}
\qquad\text{for $z\in b\O$}.\tag6.5$$
The function on the left hand side of (6.1) is holomorphic in $z$
on a neighborhood of $b\O$; so is the function on the right hand side.
Since these functions agree on $b\O$, they must be equal on a
neighborhood of $b\O$.  In fact, because $K$ and $\Lambda$ are rational,
these two functions are equal as meromorphic functions on the
neighborhood $\Cal O$ of $b\O$ on which $R(z)$ is defined.  We may
assume that $\Cal O$ is small enough that $K(z,a)$ and $\Lambda(z,a)$ have
no poles or zeroes in $\Cal O$.  Formula (6.5) now allows us to read
off the following facts.
Neither $K(z,w_k)$ nor $\Lambda(z,w_k)$ can have a  pole $z_0\in D_1$
near $\gamma_1$ because neither of these functions has a pole at the
reflected point $R(z_0)\in\O$.

Notice that (6.4) yields that
$$\Lambda(z,A_n)T(A_n)=-K(z,A_n)\overline{T(A_n)}$$
for $z\in\O$, and hence for $z$ in $D_1$.  Hence, $K(z,A_n)$ and
$\Lambda(z,A_n)$ have the same poles in $\overline{D_1}$.
Because no poles of
$K(z,w_k)$ or $\Lambda(z,w_k)$ can migrate near the boundary of $D_1$
as $w_k\to A_n$, we deduce that $P_K=P_\Lambda$.

We now take a sequence of points $w_k$ in $\O$ that tend to $A_1$ along
$\sigma$.  Formula (6.5) remains valid and, if we reason as above, we deduce
that, since
$\Lambda(z,w_k)$ has a double pole at $z=w_k$, it follows that
$K(z,w_k)$ has a double pole at the reflected point $R(w_k)$.
We now let $k\to\infty$ and use the facts that $K(z,A_1)$ and $\Lambda(z,A_1)$
have the same poles in $\overline{D_1}$ and that one of those
poles is a double pole at $A_1$.  We deduce that $P_K=P_\Lambda+2$.
This relationship contradicts the one we obtained by letting $w_k$
tend to $A_n$, and we conclude that $\O$ cannot be multiply
connected.  The proof is complete.
\enddemo

We shall now explain how to relax the smoothness assumption
that the boundary of $\O$ be $C^\infty$ smooth in Theorem~6.1.  The conclusion
about the Szeg\H o kernel holds if the boundary is only assumed to be $C^2$
smooth because, in this setting, the functions $S_a(z)$ and $L_a(z)$ extend
continuously to the boundary and $T(z)$ is continuous on $b\O$.  All of the
arguments carry through as before.

We next show that the conclusion about the Bergman kernel in Theorem~6.1
holds if the domain $\O$ is only assumed to be finitely connected and such
that no boundary component is a point.  Since the Bergman kernel is related
to the Green's function by
$K(z,w)=(-2/\pi)\frac{\dee^2}{\dee z\dee\bar w}G(z,w)$, it follows
that if the Green's function is the logarithm of a rational function
of the real and imaginary parts of $z$ and $w$, then the Bergman
kernel must be rational, too.  Hence, we have the following
theorem.

\proclaim{Theorem 6.3}
Suppose $\O$ is a finitely connected domain such that no boundary
component of $\O$ is a point.
The Green's function $G(z,a)$ associated to $\O$ is the logarithm of
a real valued rational function of the four real variables given by the
real and imaginary parts of $z$ and $a$ if and only if $\O$ is simply
connected and equivalent to the disc via a rational biholomorphic
mapping.  Similarly, the Bergman kernel $K(z,w)$ associated to $\O$ is rational
if and only if $\O$ is simply connected and equivalent to the disc via a
rational biholomorphic mapping.
\endproclaim

Hence, the only finitely connected domains having Green's functions as
simple as the Green's function for the disc are the obvious ones.
(Of course, the Green's function itself can never be rational because
it has a logarithmic singularity.)

\demo{Proof of Theorem 6.3}
Suppose $\O$ is an $n$-connected domain such that no boundary
component is a point and assume that the Bergman kernel $K(z,w)$
associated to $\O$ is rational.  It is a standard result in the theory
of conformal mapping that $\O$ is biholomorphic to a bounded domain
with real analytic boundary.  Let $\Ot$ denote such a bounded $n$-connected
domain with $C^\infty$ smooth boundary whose boundary consists
of $n$ non-intersecting simple closed real analytic curves and let
$\Phi:\O\to\Ot$ denote the biholomorphic mapping.  Let $\widetilde K(z,w)$
denote the Bergman kernel associated to $\Ot$.  The transformation formula
for the Bergman kernel under biholomorphic mappings gives
$$K(z,w)=\Phi'(z)\widetilde K(\Phi(z),\Phi(w))\overline{\Phi'(w)}.\tag6.6$$

It will be convenient to operate in the extended complex plane because it
is inconvenient if the point at infinity belongs to one of the boundary
components of $\O$.  The transformation formula for the Bergman kernel
under biholomorphic maps allows us to replace $\O$ by any domain which is
the inverse image of $\O$ under a rational biholomorphic map.  By replacing
$\O$ by its inverse image under a mapping of the form $1/(z-a)$, we may
suppose that $\O$ contains the point at infinity in its interior.

Since the transformation formula for the Bergman kernels under proper
holomorphic mappings holds in the more general setting of Theorem~6.3, we
deduce, as above, that the Ahlfors mappings are rational when the Bergman
kernel is rational.  Pick a point $a\in\O$ and let $f_a(z)$ denote the
Ahlfors map associated to $a$.  Since $f_a$ is rational, and since it
is clear that $f_a$ cannot have any poles in $\Obar$, it follows that
the boundary of $\O$ consists of finitely many piecewise real analytic
curves.  Furthermore, there are at most finitely many points in the
boundary where the boundary is not a $C^\infty$ smooth curve.  The non-smooth
points in the boundary occur at boundary points where $f_a'$ vanishes.
Suppose $f_a'$ vanishes to order $m$ at a boundary point $z_0$.  The
boundary of $\O$ near $z_0$ is described by two real analytic curves
that cross at $z_0$ and make an angle of $\pi/(m+1)$.  The mapping
$\Phi:\O\to\Ot$ described above extends continuously to the boundary
of $\O$.  Let $A=\Phi(a)$, and let $F_A(z)$ denote the
Ahlfors map of $\Ot$ onto the unit disc associated to $A$.  Since Ahlfors
maps are solutions to an extremal problem of mapping the domain into
the unit disc in such a way so as to maximize the real part of the derivative
of the mapping at the associated point, it is easy to see that Ahlfors maps
are invariant under biholomorphic mappings modulo unimodular constants to
make derivatives real valued at the points of interest.  Hence, we may write
$$f_a=e^{i\theta}F_A\circ\Phi,$$
where $\theta$ is a real constant.  Since $\Ot$ has real analytic boundary,
the Ahlfors map $F_A$ extends holomorphically past the boundary and is locally
one-to-one near the boundary.  Hence, near $z_0$, we may write
$$\Phi=F_A^{-1}\circ (e^{-i\theta}f_a)$$
to see that $\Phi$ extends holomorphically past the boundary of $\O$
near $z_0$ and $\Phi'$ vanishes to order $m$ at $z_0$.  Hence $\Phi$
extends holomorphically to a neighborhood of $\Obar$ and $\Phi'$ only
vanishes at points in $\Obar$ that are corners in the boundary.
Formula (6.6) now yields that $K(z,w)$ extends holomorphically in $z$
and antiholomorphically in $w$ to a neighborhood in $\C\times\C$ of
$(\Obar\times\Obar)-\{(z,z):z\in b\O\}$.  The rest of the argument is
now a routine transcription of the proof of Theorem~6.1.  All of the
kernel identities used in the proof of Theorem~6.1 can be deduced by
pulling back the identities that are known on $\Ot$.  For example, the
fact that $|f_a(z)|=|f_z(a)|$ can be deduced by using the argument
given in the proof of Theorem~6.1 and then pulling back to $\O$ using
$\Phi$.  The Green's function $G(z,w)$ is related to the Green's
function $\widetilde G(z,w)$ on $\Ot$ via $G(z,w)=G(\Phi(z),\Phi(w))$ and
the corresponding statement for the $\Lambda$ kernels is
$\Lambda(x,w)=\Phi'(z)\widetilde\Lambda(\Phi(z),\Phi(w))\Phi'(w)$.  The
movement of zeroes and poles of $K(z,w)$ near the boundary of $\O$
can be read off from (6.6) and the known behavior of the zeroes and poles
of $\widetilde K(z,w)$ near the boundary of $\Ot$.  The kernel
$K(z,w)$ vanishes identically when $z=z_0$ is a corner in the boundary,
but this does not interfere with our work because we may choose a
curve $\sigma$ as in the proof of Theorem~6.1 that does not begin or
terminate at a corner in the boundary of $\O$.  As $w$ moves along
such a curve, the poles of $K(z,w)$ as a function of $z$ that lie in
a bounded component $\overline{D_1}$ of the complement of $\O$ cannot
approach a corner in $b\O$.  We leave it to the reader to complete the
proof.
\enddemo

\subhead 7. Complexity of the Poisson kernel
\endsubhead
I showed in \cite{4} how the Szeg\H o projection can be used to solve
the Dirichlet problem.  The method gives rise to a formula for the
Poisson kernel of a bounded $n$-connected domain $\O$ with $C^\infty$
smooth boundary which, in light of results in \S4, reveals the level
of complexity of that kernel.  We shall use the same notation for
describing $\O$ as we have set up previously, and as before, we
also select a point $a\in\O$ such that the zeroes $a_1,\dots,a_{n-1}$ of
$S(z,a)$ are all distinct and simple.  As before, let $S_a(z)=S(z,a)$ and
$L_a(z)=L(z,a)$.  The Szeg\H o projection $P$ associated to $\O$ is the
orthogonal projection of $L^2(b\O)$ onto the Hardy space $H^2(b\O)$.
The Szeg\H o kernel is the kernel for the Szeg\H o projection in the
sense that, given a function $u\in L^2(b\O)$, the projection $Pu$ is
identified with a holomorphic function $h=Pu$ defined on $\O$ whose
$L^2$ boundary values are equal to $Pu$, and
$$(Pu)(z)=\int_{w\in b\O}S(z,w)\,u(w)\ ds.$$
The Szeg\H o projection maps $C^\infty(b\O)$ into $C^\infty(\Obar)$
(see \cite{2} for proofs of these basic facts).

Recall that the set of functions $\{L(z,a_k)S(z,a)\}_{k=1}^{n-1}$
spans the same linear space as the set of functions $\{F_k'\}_{k=1}^{n-1}$.
Define an $n-1\times n-1$ matrix of periods via
$$A_{jk}=-i\int_{\gamma_j}L(z,a_k)S(z,a)\,dz,\tag7.1$$
for $j=1,\dots,n-1$.  Because the matrix of periods of $F_k'$ is
non-singular, so is $[A_{jk}]$.  The following theorem was proved in
\cite{4}.

\proclaim{Theorem 7.1}
Given $\phi\in C^\infty(b\O)$, let $c_j$ solve the linear system
$$\sum_{j=1}^{n-1}A_{jk}c_j=P(S_a\phi)(a_k),\qquad k=1,\dots,n-1.$$
The harmonic extension $\Cal E\phi$ of $\phi$ to $\O$ is given by
$$\Cal E\phi=h+\overline{H}+\sum_{j=1}^{n-1}c_j\omega_j,$$
where, if we let $\psi=\phi-\sum_{j=1}^{n-1}c_j\omega_j$, then
$$h=\frac{P(S_a\psi)}{S_a}$$
and
$$H=\frac{P(L_a\overline{\psi}\,)}{L_a}.$$
The functions $h$ and $H$ are in $A^\infty(\O)$.
\endproclaim

This theorem allows the Poisson kernel to be written down in terms of
the Szeg\H o and Garabedian kernels.  Let $[B_{jk}]$ denote the inverse of
$[A_{jk}]$ so that $c_j=\sum_{k=1}^{n-1}B_{jk}P(S_a\phi)(a_k)$, i.e., so that
$$c_j=\int_{w\in b\O}
\left(\sum_{k=1}^{n-1}B_{jk}S(a_k,w)S(w,a)\right)
\phi(w)\ ds.$$
The formulas for $h$ and $H$ can be written
$$h(z)=\int_{w\in b\O} \frac{S(z,w)S(w,a)}{S(z,a)}\psi(w)\ ds,$$
and
$$H(z)=\int_{w\in b\O} \frac{S(z,w)L(w,a)}{L(z,a)}\ \overline{\psi(w)}\ ds.$$
Finally, when all these formulas are collected in one sum, we see that
the Poisson extension $\Cal E u$ of $u$ to $\O$ is given by an integral
$$(\Cal E u)(z)=\int_{w\in b\O} p(z,w)\,u(w)\ ds,$$
where $p(z,w)$ is the Poisson kernel and is given by
$$\gather
p(z,w)=
\frac{S(z,w)S(w,a)}{S(z,a)}+
\frac{\overline{S(z,w)L(w,a)}}{\overline{L(z,a)}} \\
-\sum_{j,k=1}^{n-1}\left(B_{jk}S(a_k,w)S(w,a)
\int_{\zeta\in\gamma_j} \frac{S(z,\zeta)S(\zeta,a)}{S(z,a)}\ ds\right) \\
-\sum_{j,k=1}^{n-1}\left(\overline{B_{jk}S(a_k,w)S(w,a)}
\int_{\zeta\in\gamma_j}
\frac{\overline{S(z,\zeta)L(\zeta,a)}}{\overline{L(z,a)}}\ ds\right) \\
+\sum_{j=1}^{n-1}
\omega_j(z) \left(\sum_{k=1}^{n-1}B_{jk}S(a_k,w)S(w,a)\right).
\endgather$$

A disappointing feature of this formula for the Poisson kernel is the
appearance of the term $\omega_j(z)$.  This function is closely tied
to the Szeg\H o kernel, but it is not as easily computed as the other
terms in the sum.  I gave a method to compute $F_j'$ in \cite{4, page~12}.
The function $\omega_j$ can be gotten from $F_j'$ via the identity
$$\omega_j(z)=\frac{1}{2\pi i}\iint_{w\in\O}\frac{F_j'(w)}{w-z}
\ dw\wedge d\bar w.$$
This is the one point in this paper where I have not been able to
obviate the need to compute an integral with respect to area
measure.

\subhead 8. Formulas for other kernels \endsubhead
The formulas for the Szeg\H o kernel and Bergman kernel are the most
interesting results of this paper.  Similar formulas may be deduced
for the Garabedian kernel $L(z,w)$ and the kernel $\Lambda(z,w)$.
We close this paper by writing these formulas down.

Let $z\in\O$ and $w\in b\O$, and consider formula (3.2).
Using identity (2.1) and the fact that $\overline{f_a}=1/f_a$ on $b\O$,
we obtain
$$L(z,w)=\frac{f_a(w)}{f_a(z)-f_a(w)}\left(c_0 S(z,a)L(w,a)+
\sum_{i,j=1}^{n-1} c_{ij}S(z,a_i)L(w,a_j)\right).$$
Since both sides of this identity are holomorphic in $z$ and $w$, this
identity holds for $z,w\in\O$, $z\ne w$.  Note that the constants
$c_0$ and $c_{ij}$ are the same as the constants in (3.2).

Similarly, combining (6.4), (5.1), and (2.1) yields the identity
$$\Lambda(z,w)=4\pi L(z,w)^2
-\sum_{i,j=1}^{n-1}\lambda_{ij}S(z,a_i)L(z,a)L(w,a_j)S(w,a),$$
where the coefficients $\lambda_{ij}$ are the same as those appearing
in (5.1), and $z\in\O$ and $w\in b\O$.  Again, since both sides of this
identity are holomorphic in $z$ and $w$, the identity holds for $z,w\in\O$,
$z\ne w$.

\enddocument